\documentclass[preprint,12pt]{elsarticle}
\usepackage{amssymb}
\usepackage{amsmath}

\newtheorem{thm}{Theorem}[section]
\newtheorem{lem}[thm]{Lemma}
\newtheorem{cor}[thm]{Corollary}
\newdefinition{defn}[thm]{Definition}
\newdefinition{rem}[thm]{Remark}
\newproof{pf}{Proof}
\newproof{pot1}{Proof of Theorem \ref{slj}}
\newproof{pot2}{Proof of Theorem \ref{czx}}
\newproof{pot3}{Proof of Theorem \ref{wqj}}

\journal{arXiv}
\begin{document}
\begin{frontmatter}

\title{Existence and Multiplicity of Nontrivial Weak Solutions for Kirchhoff-type Fractional $p$-Laplacian Equation}

\author{Taiyong Chen,\ \ Wenbin Liu
\footnote{Corresponding author.\\
{\it Telephone number:} (86-516) 83591530. {\it Fax number:} (86-516) 83591591.
{\it E-mail addresses:} taiyongchen@cumt.edu.cn (T. Chen), wblium@163.com (W. Liu), jinhua197927@163.com (H. Jin).}
,\ \ Hua Jin}

\address{Department of Mathematics, China University of Mining and Technology, Xuzhou 221116, PR China}

\begin{abstract}
We discuss the Kirchhoff-type $p$-Laplacian Dirichlet problem containing the left and right fractional derivative operators. By using the mountain pass theorem and the genus properties in critical point theory, we get some new results on the existence and multiplicity of nontrivial weak solutions for such Dirichlet problem.
\end{abstract}


\end{frontmatter}

\section{Introduction}
\label{sec1}
In this paper, we are concerned with the existence and multiplicity of nontrivial weak solutions for the Kirchhoff-type fractional Dirichlet problem with $p$-Laplacian of the form
\begin{eqnarray}
\label{kbvp}
\left\{
\begin{array}{ll}
\left(a+b\int_0^T|{_0D_t^\alpha}u(t)|^pdt\right)^{p-1}{_tD_T^\alpha}\phi_p({_0D_t^\alpha}u(t))=f(t,u(t)),\ \ t\in(0,T),\\
u(0)=u(T)=0,
\end{array}
\right.
\end{eqnarray}
where $a,b>0,\ p>1$ are constants, $_0D_t^\alpha$ and $_tD_T^\alpha$ are the left and right Riemann-Liouville fractional derivatives of order $\alpha\in(1/p,1]$ respectively, $\phi_p:\mathbb{R}\rightarrow\mathbb{R}$ is the $p$-Laplacian defined by
\begin{eqnarray*}
\phi_p(s)=|s|^{p-2}s\ (s\neq0),\ \ \phi_p(0)=0,
\end{eqnarray*}
and $f\in C([0,T]\times\mathbb{R},\mathbb{R})$.

The Kirchhoff equation introduced by Kirchhoff (\cite{gy18}) is an extension of the D'Alembert wave equation for free vibrations of elastic strings, which considers the changes in length of the string produced by transverse vibrations. Moreover, due to the fractional derivatives provide an excellent tool to describe the memory and hereditary properties of various processes and materials (\cite{dkfa,3,hilf,kjfx,main}), the fractional order models are more adequate than the integer order models in some real world problems. In addition, the $p$-Laplacian presented by Leibenson (\cite{lsl}) often occurs in non-Newtonian fluid theory, nonlinear elastic mechanics and so forth.

Note that, when $a=1$, $b=0$ and $p=2$, the left-hand side of equation of BVP (\ref{kbvp}), which is nonlinear and nonlocal, reduces to the linear differential operator $_tD_T^\alpha{_0D_t^\alpha}$, and further reduces to the local second-order differential operator $-d^2/dt^2$ when $\alpha=1$.

In the past few decades, many important results on the fractional boundary value problems (BVPs for short) (\cite{21,zbhl,23,e1,25,wj}) and the Kirchhoff equations (\cite{gy2,gy4,gy6,gy11,gy19,gy25}) have been obtained. Moreover, since the applications in physical phenomena exhibiting anomalous diffusion, the models containing left and right fractional differential operators are recently obtaining more attention (\cite{jbe3,jcre,jerw,12,fjy,zzr}).

Motivated by the above works, the purpose of this paper is to study the existence and multiplicity of nontrivial weak solutions for BVP (\ref{kbvp}). More precisely, we prove that BVP (\ref{kbvp}) possesses at least one nontrivial weak solution when $f(t,x)$ is $(p^2-1)$-superlinear or $(p^2-1)$-sublinear in $x$ at infinity, and possesses infinitely many nontrivial weak solutions when $f(t,x)$ is $(p^2-1)$-sublinear in $x$ at infinity. The main ingredients used here are the mountain pass theorem and the genus properties in critical point theory. Note that, since the Kirchhoff-type $p$-Laplacian is a nonlinear operator, it is usually difficult to verify the Palais-Smale condition ((PS)-condition for short).

In order to state our main results, we make the following assumptions on the nonlinearity $f$.

$(H_{11})$ There exist two constants $\mu>p^2,\ R>0$ such that
\begin{eqnarray*}
0<\mu F(t,x)\leq xf(t,x),\ \ \forall t\in[0,T],\ x\in\mathbb{R}\ \mbox{with}\ |x|\geq R,
\end{eqnarray*}
where $F(t,x)=\int_0^xf(t,s)ds$.

$(H_{12})$ $f(t,x)=o(|x|^{p-1})$ as $|x|\rightarrow0$ uniformly for $\forall t\in[0,T]$.

$(H_{21})$ There exist a constant $1<r_1<p^2$ and a function $d\in L^1([0,T],\mathbb{R}^+)$ such that
\begin{eqnarray*}
|f(t,x)|\leq r_1a(t)|x|^{r_1-1},\ \ \forall (t,x)\in[0,T]\times\mathbb{R}.
\end{eqnarray*}

$(H_{22})$ There exist an open interval $\mathbb{I}\subset[0,T]$ and three constants $\eta,\delta>0$, $1<r_2<p$ such that
\begin{eqnarray*}
F(t,x)\geq\eta|x|^{r_2},\ \ \forall (t,x)\in\mathbb{I}\times[-\delta,\delta].
\end{eqnarray*}

$(H_{23})$ $f(t,x)=-f(t,-x),\ \ \forall (t,x)\in[0,T]\times\mathbb{R}$.\\
Note that the assumption $(H_{11})$ implies
\begin{eqnarray*}
F(t,x)\geq c_1|x|^\mu-c_2,\ \ \forall (t,x)\in[0,T]\times\mathbb{R},
\end{eqnarray*}
where $c_1,c_2>0$ are two constants, that is, $f(t,x)$ is $(p^2-1)$-superlinear in $x$ at infinity. However $f(t,x)$ satisfying $(H_{21})$ is $(p^2-1)$-sublinear in $x$ at infinity. In addition $F(t,x)$ satisfying $(H_{11})$ is supposed to be nonnegative. However $F(t,x)$ satisfying $(H_{21})$-$(H_{23})$ can change its sign.

We are now in a position to state our main results.

\begin{thm}
\label{slj}
Let $(H_{11})$ and $(H_{12})$ be satisfied. Then BVP (\ref{kbvp}) possesses at least one nontrivial weak solution.
\end{thm}

\begin{thm}
\label{czx}
Let $(H_{21})$ and $(H_{22})$ be satisfied. Then BVP (\ref{kbvp}) possesses at least one nontrivial weak solution.
\end{thm}

\begin{thm}
\label{wqj}
Let $(H_{21})$-$(H_{23})$ be satisfied. Then BVP (\ref{kbvp}) possesses infinitely many nontrivial weak solutions.
\end{thm}

The remaining part of this paper is organized as follows. Some preliminary results are presented in Section \ref{sec2}. Section \ref{sec3} and \ref{sec4} are devoted to carry out the proof of Theorem \ref{slj}-\ref{wqj}.

\section{Preliminaries}
\label{sec2}

\subsection{Fractional Sobolev space}

By means of the critical point theory, to get the existence of solutions for BVPs, a suitable underlying function space is necessary. In this subsection, we introduce a fractional Sobolev space and some properties of this space (\cite{fjy}). Moreover we present some basic definitions and notations of the fractional calculus (\cite{15,18}).

\begin{defn}%[\cite{15}]
\label{defn2.1}
For $\gamma>0$, the left and right Riemann-Liouville fractional integrals of order $\gamma$ of a function $u:[a,b]\rightarrow\mathbb{R}$ are given by
\begin{eqnarray*}
&&_aI_{t}^\gamma u(t)=\frac{1}{\Gamma(\gamma)}\int_a^t(t-s)^{\gamma -1}u(s)ds,\\
&&_tI_{b}^\gamma u(t)=\frac{1}{\Gamma(\gamma)}\int_t^b(s-t)^{\gamma -1}u(s)ds,
\end{eqnarray*}
provided that the right-hand side integrals are pointwise defined on $[a,b]$, where $\Gamma(\cdot)$ is the Gamma function.
\end{defn}

\begin{defn}%[\cite{15}]
\label{defn2.2}
For $n-1\leq\gamma<n\ (n\in\mathbb{N})$, the left and right Riemann-Liouville fractional derivatives of order $\gamma$ of a function $u:[a,b]\rightarrow\mathbb{R}$ are given by
\begin{eqnarray*}
&&_aD_{t}^\gamma u(t)=\frac{d^n}{dt^n}{_a}I_{t}^{n-\gamma} u(t),\\
&&_tD_{b}^\gamma u(t)=(-1)^n\frac{d^n}{dt^n}{_t}I_{b}^{n-\gamma} u(t).
\end{eqnarray*}
\end{defn}

\begin{rem}
\label{zxzj}
When $\gamma=1$, one can obtain from Definition \ref{defn2.1} and \ref{defn2.2} that
\begin{eqnarray*}
_aD_{t}^1 u(t)=u'(t),\ \ _tD_{b}^1 u(t)=-u'(t),
\end{eqnarray*}
where $u'$ is the usual first-order derivative of $u$.
\end{rem}

\begin{defn}%[\cite{cjh}]
\label{defn3.1}
For $0<\alpha\leq1$ and $1<p<\infty$, the fractional derivative space $E{_0^{\alpha,p}}$ is defined by the closure of $C_0^\infty([0,T],\mathbb{R})$ with respect to the following norm
\begin{eqnarray*}
\|u\|_{E^{\alpha,p}}=(\|u\|_{L^p}^p+\|{_0D_t^\alpha}u\|_{L^p}^p)^{\frac{1}{p}},
\end{eqnarray*}
where $\|u\|_{L^p}=\left(\int_0^T|u(t)|^pdt\right)^{1/p}$ is the norm of $L^p([0,T],\mathbb{R})$.
\end{defn}

\begin{rem}
It is obvious that, for $u\in E{_0^{\alpha,p}}$, one has
\begin{eqnarray*}
u,{_0D_t^\alpha}u\in L^p([0,T],\mathbb{R}),\ \ u(0)=u(T)=0.
\end{eqnarray*}
\end{rem}

\begin{lem}[see \cite{fjy}]
\label{lem1}
Let $0<\alpha\leq1$ and $1<p<\infty$. The fractional derivative space $E{_0^{\alpha,p}}$ is a reflexive and separable Banach space.
\end{lem}

\begin{lem}[see \cite{fjy}]
Let $0<\alpha\leq1$ and $1<p<\infty$. For $u\in E_0^{\alpha,p}$, one has
\begin{eqnarray}
\label{cp}
\|u\|_{L^p}\leq C_p\|{_0D_t^\alpha}u\|_{L^p},
\end{eqnarray}
where
\begin{eqnarray*}
C_p=\frac{T^\alpha}{\Gamma(\alpha+1)}>0
\end{eqnarray*}
is a constant. Moreover, if $\alpha>1/p$, then
\begin{eqnarray}
\label{cwq}
\|u\|_\infty\leq C_\infty\|{_0D_t^\alpha}u\|_{L^p},
\end{eqnarray}
where $\|u\|_\infty=\max_{t\in[0,T]}|u(t)|$ is the norm of $C([0,T],\mathbb{R})$ and
\begin{eqnarray*}
C_\infty=\frac{T^{\alpha-\frac{1}p}}{\Gamma(\alpha)(\alpha q-q+1)^{\frac{1}q}}>0,\ \
q=\frac{p}{p-1}>1
\end{eqnarray*}
are two constants.
\end{lem}

\begin{rem}
By (\ref{cp}), we can consider the space $E_0^{\alpha,p}$ with norm 
\begin{eqnarray}
\label{fsdj}
\|u\|_{E^{\alpha,p}}=\|{_0D_t^\alpha}u\|_{L^p}
\end{eqnarray}
in what follows.
\end{rem}

\begin{lem}[see \cite{fjy}]
\label{thm3.2}
Let $1/p<\alpha\leq1$ and $1<p<\infty$. The imbedding of $E_0^{\alpha,p}$ in $C([0,T],\mathbb{R})$ is compact.
\end{lem}

\subsection{Critical point theory}

Now we present some necessary definitions and theorems of the critical point theory (\cite{dl1,dl2}). Let $X$ be a real Banach space, $I\in C^1(X,\mathbb{R})$ which means that $I$ is a continuously Fr\'{e}chet differentiable functional defined on $X$. Moreover let $B_\rho(0)$ be the open ball in $X$ with the radius $\rho$ and centered at $0$, and $\partial B_\rho(0)$ denote its boundary.

\begin{defn}
Let $I\in C^1(X,\mathbb{R})$. If any sequence $\{u_k\}\subset X$ for which $\{I(u_k)\}$ is bounded and $I'(u_k)\rightarrow0$ as $k\rightarrow\infty$ possesses a convergent subsequence in $X$, then we say that $I$ satisfies the (PS)-condition.
\end{defn}

\begin{lem}[\cite{dl2}]
\label{sldl}
Let $X$ be a real Banach space and $I\in C^1(X,\mathbb{R})$ satisfying the (PS)-condition. Suppose that $I(0)=0$ and\\
($C_1$) there are constants $\rho,\sigma>0$ such that $I|_{\partial B_\rho(0)}\geq\sigma$,\\
($C_2$) there is an $e\in X\backslash\overline{B_\rho(0)}$ such that $I(e)\leq0$.\\
Then $I$ possesses a critical value $c\geq\sigma$. Moreover $c$ can be characterized as
\begin{eqnarray*}
c=\inf_{\gamma\in\Gamma}\max_{s\in[0,1]}I(\gamma(s)),
\end{eqnarray*}
where
\begin{eqnarray*}
\Gamma=\{\gamma\in C([0,1],X)|\gamma(0)=0,\ \gamma(1)=e\}.
\end{eqnarray*}
\end{lem}

\begin{lem}[see \cite{dl1}]
\label{dll1}
Let $X$ be a real Banach space and $I\in C^1(X,\mathbb{R})$ satisfies the (PS)-condition. If $I$ is bounded from blow, then $c=\inf_X I$ is a critical value of $I$.
\end{lem}

For the sake of finding infinitely many critical points of $I$, the following {\it genus} properties are needed in our paper. Let
\begin{eqnarray*}
&&\Sigma=\{A\subset X-\{0\}|A\mbox{ is closed in }X\mbox{ and symmetric with respect to }0\},\\
&&K_c=\{u\in X|I(u)=c,\ I'(u)=0\},\ \ I^c=\{u\in X|I(u)\leq c\}.
\end{eqnarray*}

\begin{defn}
For $A\in\Sigma$, we say the genus of $A$ is $n$ denoted by $\gamma(A)=n$ if there is an odd map $G\in C(A,\mathbb{R}^n\backslash\{0\})$ and $n$ is the smallest integer with this property.
\end{defn}

\begin{lem}[see \cite{dl2}]
\label{dll2}
Let $I$ be an even $C^1$ functional on $X$ and satisfy the (PS)-condition. For any $n\in\mathbb{N}$, set
\begin{eqnarray*}
\Sigma_n=\{A\in\Sigma|\gamma(A)\geq n\},\ \ c_n=\inf_{A\in\Sigma_n}\sup_{u\in A}I(u).
\end{eqnarray*}
(i) If $\Sigma_n\neq\emptyset$ and $c_n\in\mathbb{R}$, then $c_n$ is a critical value of $I$.\\
(ii) If there exists $l\in\mathbb{N}$ such that $c_n=c_{n+1}=\cdots=c_{n+l}=c\in\mathbb{R}$, and $c\neq I(0)$, then $\gamma(K_c)\geq l+1$.
\end{lem}

\begin{rem}
From Remark 7.3 in \cite{dl2}, we know that if $K_c\in\Sigma$ and $\gamma(K_c)>1$, then $K_c$ contains infinitely many distinct points, that is, $I$ has infinitely many distinct critical points in $X$.
\end{rem}

\section{Proof of Theorem \ref{slj}}
\label{sec3}

In this section, we discuss the existence of nontrivial weak solutions of BVP (\ref{kbvp}) when the nonlinearity $f(t,x)$ is $(p^2-1)$-superlinear in $x$ at infinity. To this end, we are going to set up the corresponding variational framework of BVP (\ref{kbvp}).

Define the functional $I:E_0^{\alpha,p}\rightarrow\mathbb{R}$ by
\begin{align}
\label{fh}
I(u)
&=\frac{1}{bp^2}\left(a+b\int_0^T|{_0D_t^\alpha}u(t)|^pdt\right)^p-\int_0^TF(t,u(t))dt
-\frac{a^p}{bp^2}\nonumber\\
&=\frac{1}{bp^2}(a+b\|u\|_{E^{\alpha,p}}^p)^p-\int_0^TF(t,u(t))dt-\frac{a^p}{bp^2}.
\end{align}
In this paper, by the weak solutions of BVP (\ref{kbvp}) we mean the critical points of the associated energy functional $I$. It is easy to verify from (\ref{cwq}), (\ref{fsdj}) and $f\in C([0,T]\times\mathbb{R},\mathbb{R})$ that the functional $I$ is well defined on $E_0^{\alpha,p}$ and is a continuously Fr\'{e}chet differentiable functional, that is, $I\in C^1(E_0^{\alpha,p},\mathbb{R})$. Furthermore we have
\begin{align}
\label{fhds}
\langle I'(u),v\rangle
&=(a+b\|u\|_{E^{\alpha,p}}^p)^{p-1}
\int_0^T\phi_p({_0D_t^\alpha}u(t)){_0D_t^\alpha}v(t)dt\nonumber\\
&\ \ \ \ -\int_0^Tf(t,u(t))v(t)dt,\ \ \forall u,v\in E_0^{\alpha,p},
\end{align}
which yields
\begin{eqnarray}
\label{fhds2}
\langle I'(u),u\rangle
=(a+b\|u\|_{E^{\alpha,p}}^p)^{p-1}\|u\|_{E^{\alpha,p}}^p-\int_0^Tf(t,u(t))u(t)dt.
\end{eqnarray}

In the following, for simplicity, let
\begin{eqnarray*}
M_{pk}=a+b\|u_k\|_{E^{\alpha,p}}^p,\ \ M_p=a+b\|u\|_{E^{\alpha,p}}^p.
\end{eqnarray*}

\begin{lem}
\label{pstj}
Assume that $(H_{11})$ holds. Then $I$ satisfies the (PS)-condition in $E_0^{\alpha,p}$.
\end{lem}

\begin{pf}
Let $\{u_k\}\subset E_0^{\alpha,p}$ be a sequence such that
\begin{eqnarray*}
|I(u_k)|\leq K,\ \ I'(u_k)\rightarrow0\ \mbox{as}\ {k\rightarrow\infty},
\end{eqnarray*}
where $K>0$ is a constant. We first prove that $\{u_k\}$ is bounded in $E_0^{\alpha,p}$. From the continuity of $\mu F(t,x)-xf(t,x)$ and $(H_{11})$, we obtain that there exists a constant $c>0$ such that
\begin{eqnarray*}
F(t,x)\leq\frac{1}{\mu}xf(t,x)+c,\ \ \forall (t,x)\in[0,T]\times\mathbb{R}.
\end{eqnarray*}
Thus, by (\ref{fh}) and (\ref{fhds2}), we have
\begin{align*}
K&\geq I(u_k)\\
&\geq\frac{1}{bp^2}M_{pk}^p
-\frac{1}{\mu}\int_0^Tf(t,u_k(t))u_k(t)dt-cT-\frac{a^p}{bp^2}\\
&=\frac{1}{bp^2}M_{pk}^p
-\frac{1}{\mu}M_{pk}^{p-1}\|u_k\|_{E^{\alpha,p}}^p+\frac{1}{\mu}\langle I'(u_k),u_k\rangle-cT-\frac{a^p}{bp^2}\\
&\geq M_{pk}^{p-1}\left(\left(\frac{1}{p^2}-\frac{1}{\mu}\right)
\|u_k\|_{E^{\alpha,p}}^p+\frac{a}{bp^2}\right)\\
&\ \ \ \ -\frac{1}{\mu}\|I'(u_k)\|_{(E_0^{\alpha,p})^*}\|u_k\|_{E^{\alpha,p}}-cT-\frac{a^p}{bp^2},
\end{align*}
which together with $I'(u_k)\rightarrow0\ \mbox{as}\ {k\rightarrow\infty}$ yields
\begin{align*}
K&\geq M_{pk}^{p-1}\left(\left(\frac{1}{p^2}-\frac{1}{\mu}\right)
\|u_k\|_{E^{\alpha,p}}^p+\frac{a}{bp^2}\right)-\|u_k\|_{E^{\alpha,p}}\\
&\ \ \ \ -cT-\frac{a^p}{bp^2}\ \ \mbox{as}\ k\rightarrow\infty.
\end{align*}
Then it follows from $\mu>p^2$ that $\{u_k\}$ is bounded in $E_0^{\alpha,p}$.

Since $E_0^{\alpha,p}$ is a reflexive Banach space (see Lemma \ref{lem1}), going if necessary to a subsequence, we can assume $u_k\rightharpoonup u$ in $E_0^{\alpha,p}$. Hence, from $I'(u_k)\rightarrow0\ \mbox{as}\ {k\rightarrow\infty}$ and the definition of weak convergence, we have
\begin{align}
\label{jl3}
&\langle I'(u_k)-I'(u),u_k-u\rangle \nonumber\\
&=\langle I'(u_k),u_k-u\rangle -\langle I'(u),u_k-u\rangle \nonumber\\
&\leq\|I'(u_k)\|_{(E_0^{\alpha,p})^*}\|u_k-u\|_{E^{\alpha,p}}-\langle I'(u),u_k-u\rangle \nonumber\\
&\rightarrow0\ \ \mbox{as}\ k\rightarrow\infty.
\end{align}
In addition, we obtain from (\ref{cwq}), (\ref{fsdj}) and Lemma \ref{thm3.2} that $\{u_k\}$ is bounded in $C([0,T],\mathbb{R})$ and $\|u_k-u\|_\infty\rightarrow0$ as $k\rightarrow\infty$. Thus there exists a constant $c_1>0$ such that
\begin{eqnarray*}
|f(t,u_k(t))-f(t,u(t))|\leq c_1,\ \ \forall t\in[0,T],
\end{eqnarray*}
which yields
\begin{align}
\label{jl4}
\left|\int_0^T(f(t,u_k(t))-f(t,u(t)))(u_k(t)-u(t))dt\right|
&\leq c_1T\|u_k-u\|_\infty\nonumber\\
&\rightarrow0\ \ \mbox{as}\ k\rightarrow\infty.
\end{align}
Moreover, by the boundedness of $\{u_k\}$ in $E_0^{\alpha,p}$, one has
\begin{align}
\label{528.1}
&(M_{pk}^{p-1}-M_p^{p-1})
\int_0^T\phi_p({_0D_t^\alpha}u(t))({_0D_t^\alpha}u_k(t)-{_0D_t^\alpha}u(t))dt\nonumber\\
&=(M_{pk}^{p-1}-M_p^{p-1})\langle I_1'(u),u_k-u\rangle\nonumber\\
&\rightarrow0\ \ \mbox{as}\ k\rightarrow\infty,
\end{align}
where $I_1'$ is the Fr\'{e}chet derivative of $I_1:E_0^{\alpha,p}\rightarrow\mathbb{R}$ defined by
\begin{eqnarray*}
I_1(u)=\frac{1}{p}\int_0^T|{_0D_t^\alpha}u(t)|^pdt.
\end{eqnarray*}
From (\ref{fhds}), we have
\begin{align*}
&\langle I'(u_k)-I'(u),u_k-u\rangle+\int_0^T(f(t,u_k(t))-f(t,u(t)))(u_k(t)-u(t))dt\\
&=M_{pk}^{p-1}\int_0^T\phi_p({_0D_t^\alpha}u_k(t))({_0D_t^\alpha}u_k(t)-{_0D_t^\alpha}u(t))dt\\
&\ \ \ \ -M_p^{p-1}\int_0^T\phi_p({_0D_t^\alpha}u(t))({_0D_t^\alpha}u_k(t)-{_0D_t^\alpha}u(t))dt\\
&=M_{pk}^{p-1}\int_0^T(\phi_p({_0D_t^\alpha}u_k(t))-\phi_p({_0D_t^\alpha}u(t)))
({_0D_t^\alpha}u_k(t)-{_0D_t^\alpha}u(t))dt\\
&\ \ \ \ +(M_{pk}^{p-1}-M_p^{p-1})\int_0^T\phi_p({_0D_t^\alpha}u(t))
({_0D_t^\alpha}u_k(t)-{_0D_t^\alpha}u(t))dt,
\end{align*}
which together with (\ref{jl3})-(\ref{528.1}) yields
\begin{eqnarray}
\label{jl8}
\int_0^T(\phi_p({_0D_t^\alpha}u_k(t))-\phi_p({_0D_t^\alpha}u(t)))
({_0D_t^\alpha}u_k(t)-{_0D_t^\alpha}u(t))dt
\rightarrow0
\end{eqnarray}
as $k\rightarrow\infty$.

Following (2.10) in \cite{js}, there exist two constants $c_2,c_3>0$ such that
\begin{align}
\label{jl5}
&\int_0^T(\phi_p({_0D_t^\alpha}u_k(t))-\phi_p({_0D_t^\alpha}u(t)))
({_0D_t^\alpha}u_k(t)-{_0D_t^\alpha}u(t))dt\nonumber\\
&\geq
\left\{
\begin{array}{ll}
c_2\int_0^T|{_0D_t^\alpha}u_k(t)-{_0D_t^\alpha}u(t)|^pdt,\ \ p\geq2,\\
c_3\int_0^T\frac{|{_0D_t^\alpha}u_k(t)-{_0D_t^\alpha}u(t)|^2}
{(|{_0D_t^\alpha}u_k(t)|+|{_0D_t^\alpha}u(t)|)^{2-p}}dt,\ \ 1<p<2.
\end{array}
\right.
\end{align}
When $1<p<2$, based on the H\"older inequality, we get
\begin{align*}
&\int_0^T|{_0D_t^\alpha}u_k(t)-{_0D_t^\alpha}u(t)|^pdt\\
&\leq\left(\int_0^T\frac{|{_0D_t^\alpha}u_k(t)-{_0D_t^\alpha}u(t)|^2}
{(|{_0D_t^\alpha}u_k(t)|+|{_0D_t^\alpha}u(t)|)^{2-p}}dt\right)^{\frac{p}{2}}\\
&\ \ \ \ \cdot\left(\int_0^T(|{_0D_t^\alpha}u_k(t)|+|{_0D_t^\alpha}u(t)|)^pdt\right)
^{\frac{2-p}{2}}\\
&\leq c_4(\|u_k\|_{E^{\alpha,p}}^p+\|u\|_{E^{\alpha,p}}^p)^{\frac{2-p}{2}}
\left(\int_0^T\frac{|{_0D_t^\alpha}u_k(t)-{_0D_t^\alpha}u(t)|^2}
{(|{_0D_t^\alpha}u_k(t)|+|{_0D_t^\alpha}u(t)|)^{2-p}}dt\right)^{\frac{p}{2}},
\end{align*}
where $c_4=2^{(p-1)(2-p)/2}>0$ is a constant, which together with (\ref{jl5}) implies
\begin{align}
\label{jl6}
&\int_0^T(\phi_p({_0D_t^\alpha}u_k(t))-\phi_p({_0D_t^\alpha}u(t)))
({_0D_t^\alpha}u_k(t)-{_0D_t^\alpha}u(t))dt\nonumber\\
&\geq c_3c_4^{-\frac{2}{p}}(\|u_k\|_{E^{\alpha,p}}^p+\|u\|_{E^{\alpha,p}}^p)^{\frac{p-2}{p}}
\|u_k-u\|_{E^{\alpha,p}}^2,\ \ 1<p<2.
\end{align}
When $p\geq2$, by (\ref{jl5}), we have
\begin{align}
\label{jl7}
&\int_0^T(\phi_p({_0D_t^\alpha}u_k(t))-\phi_p({_0D_t^\alpha}u(t)))
({_0D_t^\alpha}u_k(t)-{_0D_t^\alpha}u(t))dt\nonumber\\
&\geq c_2\|u_k-u\|_{E^{\alpha,p}}^p,\ \ p\geq2.
\end{align}
Then it follows from (\ref{jl8}), (\ref{jl6}) and (\ref{jl7}) that
\begin{eqnarray*}
\|u_k-u\|_{E^{\alpha,p}}\rightarrow0\ \ \mbox{as}\ \ k\rightarrow\infty.
\end{eqnarray*}
Hence $I$ satisfies the (PS)-condition. $\Box$
\end{pf}

\begin{pot1}
From $(H_{12})$, there exist two constants $0<\varepsilon<1,\delta>0$ such that
\begin{eqnarray}
\label{530.1}
F(t,x)\leq\frac{(1-\varepsilon)a^{p-1}}{pC_p^p}|x|^p,
\ \ \forall t\in[0,T],\ x\in\mathbb{R}\ \mbox{with}\ |x|\leq\delta,
\end{eqnarray}
where $C_p>0$ is a constant defined in (\ref{cp}). Let $\rho=\delta/C_\infty$ and $\sigma=\varepsilon a^{p-1}\rho^p/p$, where $C_\infty>0$ is a constant defined in (\ref{cwq}). Then, by (\ref{cwq}) and (\ref{fsdj}), we have
\begin{eqnarray*}
\|u\|_\infty\leq C_\infty\|u\|_{E^{\alpha,p}}=\delta,\ \
\forall u\in E_0^{\alpha,p}\ \mbox{with}\ \|u\|_{E^{\alpha,p}}=\rho,
\end{eqnarray*}
which together with (\ref{cp}), (\ref{fsdj}), (\ref{fh}) and (\ref{530.1}) implies
\begin{align*}
I(u)&=\frac{1}{bp^2}(a+b\|u\|_{E^{\alpha,p}}^p)^p-\int_0^TF(t,u(t))dt-\frac{a^p}{bp^2}\\
&\geq\frac{a^{p-1}}{p}\|u\|_{E^{\alpha,p}}^p
-\frac{(1-\varepsilon)a^{p-1}}{pC_p^p}\int_0^T|u(t)|^pdt\\
&\geq\frac{a^{p-1}}{p}\|u\|_{E^{\alpha,p}}^p
-\frac{(1-\varepsilon)a^{p-1}}{p}\|u\|_{E^{\alpha,p}}^p\\
&=\frac{\varepsilon a^{p-1}}{p}\|u\|_{E^{\alpha,p}}^p\\
&=\sigma,\ \ \forall u\in E_0^{\alpha,p}\ \mbox{with}\ \|u\|_{E^{\alpha,p}}=\rho,
\end{align*}
which means that the condition ($C_1$) in Lemma \ref{sldl} is satisfied.

From $(H_{11})$, a simple argument using the very definition of derivative shows that there exist two constants $c_1,c_2>0$ such that
\begin{eqnarray*}
F(t,x)\geq c_1|x|^\mu-c_2,\ \ \forall (t,x)\in[0,T]\times\mathbb{R}.
\end{eqnarray*}
Then, for any $u\in E_0^{\alpha,p}\setminus\{0\},\ \xi\in\mathbb{R}^+$, we can obtain from (\ref{fh}) and $\mu>p^2$ that
\begin{align*}
I(\xi u)&=\frac{1}{bp^2}(a+b\|\xi u\|_{E^{\alpha,p}}^p)^p-\int_0^TF(t,\xi u(t))dt-\frac{a^p}{bp^2}\\
&\leq\frac{1}{bp^2}(a+b\|\xi u\|_{E^{\alpha,p}}^p)^p
-c_1\int_0^T|\xi u(t)|^\mu dt+c_2T-\frac{a^p}{bp^2}\\
&=\frac{1}{bp^2}(a+b\xi^p\|u\|_{E^{\alpha,p}}^p)^p
-c_1\xi^\mu\|u\|_{L^\mu}^\mu+c_2T-\frac{a^p}{bp^2}\\
&\rightarrow-\infty\ \ \mbox{as}\ \xi\rightarrow\infty.
\end{align*}
Thus, taking $\xi_0$ large enough and letting $e=\xi_0u$, we have $I(e)\leq0$. Hence the condition ($C_2$) in Lemma \ref{sldl} is also satisfied.

Finally, by $I(0)=0$, Lemma \ref{sldl} and \ref{pstj}, we get a critical point $u^*$ of $I$ satisfying $I(u^*)\geq\sigma>0$, and so $u^*$ is a nontrivial solution of BVP (\ref{kbvp}). $\Box$
\end{pot1}

\section{Proof of Theorem \ref{czx} and \ref{wqj}}
\label{sec4}

In this section, we discuss the existence and multiplicity of nontrivial weak solutions of BVP (\ref{kbvp}) when the nonlinearity $f(t,x)$ is $(p^2-1)$-sublinear in $x$ at infinity.

\begin{lem}
\label{xyj}
Suppose that $(H_{21})$ is satisfied. Then $I$ is bounded from below in $E_0^{\alpha,p}$.
\end{lem}

\begin{pf}
From $(H_{21})$, one has
\begin{eqnarray*}
\label{xtj}
|F(t,u)|\leq a(t)|u|^{r_1},\ \ \forall (t,u)\in[0,T]\times\mathbb{R},
\end{eqnarray*}
which together with (\ref{cwq})-(\ref{fh}) yields
\begin{align}
\label{fsdy}
I(u)
&\geq\frac{1}{bp^2}(a+b\|u\|_{E^{\alpha,p}}^p)^p-\int_0^Ta(t)|u(t)|^{r_1}dt-\frac{a^p}{bp^2}\nonumber\\
&\geq\frac{1}{bp^2}(a+b\|u\|_{E^{\alpha,p}}^p)^p-\|a\|_{L^1}\|u\|_\infty^{r_1}-\frac{a^p}{bp^2}\nonumber\\
&\geq\frac{1}{bp^2}(a+b\|u\|_{E^{\alpha,p}}^p)^p-C_\infty^{r_1}\|a\|_{L^1}\|u\|_{E^{\alpha,p}}^{r_1}
-\frac{a^p}{bp^2}.
\end{align}
Since $1<r_1<p^2$, (\ref{fsdy}) yields $I(u)\rightarrow\infty$ as $\|u\|_{E^{\alpha,p}}\rightarrow\infty$. Hence $I$ is bounded from below. $\Box$
\end{pf}

\begin{lem}
\label{ps2}
Assume that $(H_{21})$ holds. Then $I$ satisfies the (PS)-condition in $E_0^{\alpha,p}$.
\end{lem}

\begin{pf}
Let $\{u_k\}\subset E_0^{\alpha,p}$ be a sequence such that
\begin{eqnarray*}
|I(u_k)|\leq K,\ \ I'(u_k)\rightarrow0\ \mbox{as}\ {k\rightarrow\infty},
\end{eqnarray*}
where $K>0$ is a constant. Then (\ref{fsdy}) implies that $\{u_k\}$ is bounded in $E_0^{\alpha,p}$. The remainder of proof is similar to the proof of Lemma \ref{pstj}, so we omit the details. $\Box$
\end{pf}

\begin{pot2}
From Lemma \ref{dll1}, \ref{xyj} and \ref{ps2}, we obtain $c=\inf_{E_0^{\alpha,p}}I(u)$ is a critical value of $I$, that is, there exists a critical point $u^*\in E_0^{\alpha,p}$ such that $I(u^*)=c$.

Now we show $u^*\neq0$. Let $u_0\in\left(W_0^{1,2}(\mathbb{I},\mathbb{R})\cap E_0^{\alpha,p}\right)\setminus\{0\}$ and $\|u_0\|_\infty=1$, from (\ref{fh}) and $(H_{22})$, we get
\begin{align}
\label{jfl}
I(su_0)
&=\frac{1}{bp^2}(a+b\|su_0\|_{E^{\alpha,p}}^p)^p-\int_0^TF(t,su_0(t))dt-\frac{a^p}{bp^2}\nonumber\\
&=\frac{1}{bp^2}(a+bs^p\|u_0\|_{E^{\alpha,p}}^p)^p-\int_\mathbb{I}F(t,su_0(t))dt
-\frac{a^p}{bp^2}\nonumber\\
&\leq\frac{1}{bp^2}(a+bs^p\|u_0\|_{E^{\alpha,p}}^p)^p-\eta s^{r_2}\int_\mathbb{I}|u_0(t)|^{r_2}dt\nonumber\\
&\ \ \ \ -\frac{a^p}{bp^2},\ \ 0<s\leq \delta.
\end{align}
Since $1<r_2<p$, (\ref{jfl}) implies $I(su_0)<0$ for $s>0$ small enough. Then $I(u^*)=c<0$, hence $u^*$ is a nontrivial critical point of $I$, and so $u^*$ is a nontrivial solution of BVP (\ref{kbvp}). $\Box$
\end{pot2}

\begin{pot3}
From Lemma \ref{xyj} and \ref{ps2}, we obtain that $I\in C^1(E_0^{\alpha,p},\mathbb{R})$ is bounded from below and satisfies the (PS)-condition. In addition, (\ref{fh}) and $(H_{23})$ show that $I$ is even and $I(0)=0$.

Fixing $n\in\mathbb{N}$, we take $n$ disjoint open intervals $\mathbb{I}_i$ such that $\cup_{i=1}^n\mathbb{I}_i\subset\mathbb{I}$. Let $u_i\in\left(W_0^{1,2}(\mathbb{I}_i,\mathbb{R})\cap E_0^{\alpha,p}\right)\setminus\{0\}$ and $\|u_i\|_{E^{\alpha,p}}=1$, and
\begin{eqnarray*}
E_n=\mbox{span}\{u_1,u_2,\cdots,u_n\},\ \ S_n=\{u\in E_n|\|u\|_{E^{\alpha,p}}=1\}.
\end{eqnarray*}
For $u\in E_n$, there exist $\lambda_i\in\mathbb{R}$, such that
\begin{eqnarray}
\label{zks}
u(t)=\sum_{i=1}^n\lambda_iu_i(t),\ \ \forall t\in[0,T].
\end{eqnarray}
Thus we get
\begin{align}
\label{nyfs}
\|u\|_{E^{\alpha,p}}^p
&=\int_0^T|_0D_t^\alpha u(t)|^pdt
=\sum_{i=1}^n|\lambda_i|^p\int_{\mathbb{I}_i}|_0D_t^\alpha u_i(t)|^pdt\nonumber\\
&=\sum_{i=1}^n|\lambda_i|^p\int_0^T|_0D_t^\alpha u_i(t)|^pdt
=\sum_{i=1}^n|\lambda_i|^p\|u_i\|_{E^{\alpha,p}}^p\nonumber\\
&=\sum_{i=1}^n|\lambda_i|^p,\ \ \forall u\in E_n.
\end{align}
From (\ref{cwq})-(\ref{fh}), (\ref{zks}) and $(H_{22})$, for $u\in S_n$, one has
\begin{align}
\label{fnl}
I(su)
&=\frac{1}{bp^2}(a+b\|su\|_{E^{\alpha,p}}^p)^p-\int_0^TF(t,su(t))dt-\frac{a^p}{bp^2}\nonumber\\
&=\frac{1}{bp^2}(a+bs^p)^p
-\sum_{i=1}^n\int_{\mathbb{I}_i}F(t,s\lambda_iu_i(t))dt-\frac{a^p}{bp^2}\nonumber\\
&\leq\frac{1}{bp^2}(a+bs^p)^p-\eta s^{r_2}\sum_{i=1}^n|\lambda_i|^{r_2}\int_{\mathbb{I}_i}|u_i(t)|^{r_2}dt\nonumber\\
&\ \ \ \ -\frac{a^p}{bp^2}
,\ \ 0<s\leq\frac{\delta}{C_\infty\lambda^*},
\end{align}
where $\lambda^*=\max_{i\in\{1,2,\cdots,n\}}|\lambda_i|>0$ is a constant. Since $1<r_2<p$, it follows from (\ref{fnl}) that there exist constants $\epsilon,\sigma>0$ such that
\begin{eqnarray}
\label{nlz}
I(\sigma u)<-\epsilon,\ \ \forall u\in S_n.
\end{eqnarray}
Let
\begin{eqnarray*}
S_n^\sigma=\{\sigma u|u\in S_n\},\ \ \Lambda=\left\{(\lambda_1,\lambda_2,\cdots,\lambda_n)\in\mathbb{R}^n|
\sum_{i=1}^n|\lambda_i|^p<\sigma^p\right\}.
\end{eqnarray*}
Then we obtain from (\ref{nlz}) that
\begin{eqnarray*}
I(u)<-\epsilon,\ \ \forall u\in S_n^\sigma,
\end{eqnarray*}
which, together with the fact that $I$ is even and $I(0)=0$, yields
\begin{eqnarray*}
S_n^\sigma\subset I^{-\epsilon}\in\Sigma.
\end{eqnarray*}
From (\ref{nyfs}), we know that the mapping $(\lambda_1,\lambda_2,\cdots,\lambda_n)\rightarrow\sum_{i=1}^n\lambda_iu_i(t)$ from $\partial\Lambda$ to $S_n^\sigma$ is odd and homeomorphic. Hence, by some properties of the genus (see Proposition 7.5 and 7.7 in \cite{dl2}), we deduce that
\begin{eqnarray*}
\gamma(I^{-\epsilon})\geq\gamma(S_n^\sigma)=n.
\end{eqnarray*}
Thus $I^{-\epsilon}\in \Sigma_n$ and so $\Sigma_n\neq\emptyset$. Let
\begin{eqnarray*}
c_n=\inf_{A\in\Sigma_n}\sup_{u\in A}I(u).
\end{eqnarray*}
It follows from $I$ is bounded from below that $-\infty<c_n\leq-\epsilon<0$. That is, for any $n\in\mathbb{N}$, $c_n$ is a real negative number. Hence, by Lemma \ref{dll2}, $I$ admits infinitely many nontrivial critical points, and so BVP (\ref{kbvp}) possesses infinitely many nontrivial negative energy solutions. $\Box$
\end{pot3}

Obviously the following assumption $(H_{24})$ implies $(H_{21})$-$(H_{23})$ hold.

$(H_{24})$ $f(t,x)=rg(t)|x|^{r-2}x$, where $1<r<p^2$ is a constant, $g\in C([0,T],\mathbb{R})$ and there exists an open interval $\mathbb{I}\subset[0,T]$ such that $g(t)>0,\ \forall t\in\mathbb{I}$.\\
Then a corollary can be stated as follows.

\begin{cor}
Let $(H_{24})$ be satisfied. Then BVP (\ref{kbvp}) possesses infinitely many nontrivial weak solutions.
\end{cor}

%\section*{Competing interests}
%The authors declare that they have no competing interests.

%\section*{Authors' contributions}
%The authors contributed equally in this article. They read and approved the final manuscript.

\section*{Acknowledgements}
This work was supported by the National Natural Science Foundation of China (11271364) and the Nature Science Foundation of Jiangsu Province (BK20130170).

\end{document}